\RequirePackage{ifpdf}
\ifpdf 
\documentclass[pdftex]{sigma}
\else
\documentclass{sigma}
\fi

\usepackage{epic,eepicemu}

\def\a{\alpha}
\def\b{\beta}

\def\I{{\rm I}}
\def\II{{\rm II}}
\def\III{{\rm III}}
\def\IV{{\rm IV}}
\def\V{{\rm V}}

\def\({\left(}
\def\){\right)}
\def\<{\langle}
\def\>{\rangle}

\begin{document}

\allowdisplaybreaks

\renewcommand{\thefootnote}{$\star$}

\renewcommand{\PaperNumber}{033}

\FirstPageHeading

\ShortArticleName{On Quadrirational Yang--Baxter Maps}

\ArticleName{On Quadrirational Yang--Baxter Maps\footnote{This paper is a
contribution to the Proceedings of the Workshop ``Geometric Aspects of Discrete and Ultra-Discrete Integrable Systems'' (March 30 -- April 3, 2009, University of Glasgow, UK). The full collection is
available at
\href{http://www.emis.de/journals/SIGMA/GADUDIS2009.html}{http://www.emis.de/journals/SIGMA/GADUDIS2009.html}}}

\Author{V.G.~PAPAGEORGIOU~$^{\dag^1}$,  Yu.B.~SURIS~$^{\dag^2}$, A.G.~TONGAS~$^{\dag^3}$ and A.P.~VESELOV~$^{\dag^4\dag^5}$}

\AuthorNameForHeading{V.G.~Papageorgiou,  Yu.B.~Suris, A.G.~Tongas and A.P.~Veselov}

\Address{$^{\dag^1}$~Department of Mathematics, University of Patras, 26 500 Patras, Greece}
\EmailDD{\href{mailto:vassilis@math.upatras.gr}{vassilis@math.upatras.gr}}

\Address{$^{\dag^2}$~Institut f\"ur Mathematik, Technische Universit\"at Berlin,\\
\hphantom{$^{\dag^2}$}~Str. des 17. Juni 136, 10623 Berlin, Germany}
\EmailDD{\href{mailto:suris@math.tu-berlin.de}{suris@math.tu-berlin.de}}

\Address{$^{\dag^3}$~Department of Applied Mathematics, University of Crete, 714 09 Heraklion, Greece}
\EmailDD{\href{mailto:atongas@tem.uoc.gr}{atongas@tem.uoc.gr}}

\Address{$^{\dag^4}$~School of Mathematics, Loughborough University,\\
\hphantom{$^{\dag^4}$}~Loughborough, Leicestershire, LE11 3TU, UK}
\EmailDD{\href{mailto:A.P.Veselov@lboro.ac.uk}{A.P.Veselov@lboro.ac.uk}}

\Address{$^{\dag^5}$~Moscow State University, Moscow 119899, Russia}

\ArticleDates{Received November 15, 2009, in f\/inal form March 26, 2010;  Published online April 16, 2010}

\Abstract{We use the classif\/ication of the quadrirational maps given by
Adler, Bobenko and Suris to describe when such maps satisfy the
Yang--Baxter relation. We show that the corresponding maps can be
characterized by certain singularity invariance condition. This leads
to some new families of Yang--Baxter maps corresponding to the
geometric symmetries of pencils of quadrics.}

\Classification{14E07; 14H70; 37K20}

\renewcommand{\thefootnote}{\arabic{footnote}}
\setcounter{footnote}{0}

\section{Introduction}

Recently Adler, Bobenko and Suris \cite{ABS1} introduced an
important notion of {\it quadrirational maps} of
$\mathbb{CP}^1\times\mathbb{CP}^1$ into itself. They classif\/ied
them up to a natural $({\mathcal M}\ddot{o}b)^4$-action, where
${\mathcal M}\ddot{o}b$ is the M\"obius group of the projective
transformations of $\mathbb{CP}^1,$ and gave a beautiful geometric
interpretation of the results.

In this note, which can be considered as an extended comment to
\cite{ABS1},  we use the classif\/ication of the quadrirational maps
from \cite{ABS1} to describe when such maps satisfy the
Yang--Baxter relation. We show that the corresponding maps can be
characterized by certain singularity invariance condition, which
leads to the 5~families of Yang--Baxter maps found in~\cite{ABS1},
as well as to 5~additional families corresponding to the geometric
symmetries of pencils of quadrics.

We should say that a surprising relation of the quadrirationality
with the Yang--Baxter pro\-per\-ty was discovered already in
\cite{ABS1}. However it was left without proper discussion, which
may lead to misinterpretation of the main result (see the
concluding remarks in \cite{V2}). We came naturally to this point
when  new examples of quadrirational Yang--Baxter maps were
discovered in \cite{PTV}.

\section[Yang-Baxter relation and quadrirational maps]{Yang--Baxter relation and quadrirational maps}

Let $\mathbb{X}$ be any set and $R$ a map of
$\mathbb{X}\times\mathbb{X}$ into itself. If $\mathbb{X}^n$ stands
for $\mathbb{X}\times \mathbb{X}\times \dots \times \mathbb{X}$,
then let~$R^{ij}$ will denote the map of $\mathbb{X}^n$ into
itself which acts as $R$ on the $i$-th and $j$-th factors and as
the identity on the others. A map
$R(\lambda,\mu):\mathbb{X}\times\mathbb{X}\to\mathbb{X}\times\mathbb{X}$
depending on two parameters $\lambda$, $\mu$ from some parameter set
$\Lambda$ is called a {\em Yang--Baxter map} if it satisf\/ies the
{\em Yang--Baxter relation}
\begin{gather}
R^{23} (\lambda_2, \lambda_3)  R^{13} (\lambda_1,\lambda_3)
R^{12} (\lambda_1,\lambda_2) = R^{12} (\lambda_1,\lambda_2)
R^{13} (\lambda_1,\lambda_3)  R^{23} (\lambda_2,\lambda_3),
\label{eq:YBcom}
\end{gather}
regarded as an equality of maps of $\mathbb{X} \times
\mathbb{X}\times \mathbb{X}$ into itself. If in addition the
relation
\begin{gather*}
R^{21}(\lambda_2,\lambda_1) R(\lambda_1,\lambda_2) = {\rm Id}
\end{gather*}
holds, then $R(\lambda,\mu)$ is called a {\em reversible}
Yang--Baxter map.

The following proposition, which is easily checked directly,
def\/ines a natural equivalence among the parameter-dependent
Yang--Baxter maps.

\begin{proposition}\label{Pr}
Assume that there is a family of bijections
$\phi(\lambda):\mathbb{X}\to\mathbb X$ parametrized by
$\lambda\in\Lambda$. If $R(\lambda,\mu)$ satisfies the Yang--Baxter
relation \eqref{eq:YBcom}, then the same is true for
\begin{gather}
\label{equiva} \widetilde R(\lambda,\mu) =
\phi(\lambda)^{-1}\times \phi(\mu)^{-1} R(\lambda,\mu)
\phi(\lambda)\times \phi(\mu).
\end{gather}
The Yang--Baxter maps $R$, $\widetilde R$ are called equivalent.
\end{proposition}
Parameter dependent Yang--Baxter maps can be considered as
parameter independent by extending the domain $\mathbb{X}$ to
$\mathbb{X}\times \Lambda$ with identical action of the map on the
second factor. However this may change the meaning of natural
equivalence between the maps.

From now on we will be mainly interested in the case $\mathbb
X=\mathbb{CP}^1$.

Recall that a map of $\mathbb{CP}^1\times\mathbb{CP}^1$ into
itself, $R:(x,y) \mapsto (u,v)$, is called {\it quadrirational}~\cite{ABS1} if both $R$ and the so-called {\it companion map}
$\bar R: (x, v) \mapsto (u,y)$ are birational maps. All such maps
have the form
\begin{gather}\label{maps}
 R:\  u=\frac{a(y)x+b(y)}{c(y)x+d(y)},\qquad
         v=\frac{A(x)y+B(x)}{C(x)y+D(x)},
\end{gather}
where $a(y),\ldots,d(y)$ and $A(x),\ldots,D(x)$ are polynomials of
degree at most~2. There are three subclasses of such maps, denoted
as [1:1], [1:2] and [2:2], corresponding to the highest degrees of
the coef\/f\/icients of both fractions in (\ref{maps}).

In this paper we restrict ourselves by the most interesting subclass [2:2]. Adler, Bobenko
and Suris \cite{ABS1} showed that any quadrirational map from this subclass
is $({\mathcal M}\ddot{o}b)^4$-equivalent to one of the following
f\/ive maps depending on two  complex parameters $\a$, $\b$:
\begin{alignat}{4}
\label{F1}\tag{$F_\I$}
 & u= \a yP,\qquad &&
  v= \b xP,\qquad &&
  P= \frac{(1-\b)x+\b-\a+(\a-1)y}
           {\b(1-\a)x+(\a-\b)yx+\a(\b-1)y},& \\
\label{F2}\tag{$F_\II$}
&  u= \frac{y}{\a} P, \qquad &&
  v= \frac{x}{\b} P,\qquad &&
  P= \frac{\a x-\b y+\b-\a}{x-y},& \\
\label{F3}\tag{$F_\III$}
 & u= \frac{y}{\a} P,\qquad &&
  v= \frac{x}{\b} P,\qquad &&
  P= \frac{\a x-\b y}{x-y},& \\
\label{F4}\tag{$F_\IV$}
&  u= yP,\qquad & &
  v= xP,\qquad &&
  P = 1+\frac{\b-\a}{x-y},& \\
\label{F5}\tag{$F_\V$}
&  u= y+P,\qquad &&
  v= x+P,\qquad &&
  P= \frac{\a-\b}{x-y}.&
\end{alignat}

A surprising fact is that all f\/ive canonical representative maps
satisfy the Yang--Baxter relation
\begin{gather*}
R^{23} (\a_2,\a_3)  R^{13} (\a_1,\a_3)  R^{12} (\a_1,\a_2) =
R^{12} (\a_1,\a_2)  R^{13} (\a_1,\a_3)  R^{23} (\a_2,\a_3).
\end{gather*}
However, it should be noted that not all quadrirational maps
satisfy the Yang--Baxter relation (as one might conclude from~\cite{ABS1}), since the action of $({\mathcal M}\ddot{o}b)^4$, in
general, destroys the Yang--Baxter property. For example, if we change in $F_\V$
the variables $x$, $y$, $u$, $v$ to $-x$, $-y$, $u$, $v$, we come to the map
\[
u= -(y+P),\qquad v= -(x+P),\qquad P= \frac{\a-\b}{x-y},
\]
which does not satisfy the Yang--Baxter property (see \cite{V2}).

\section[Yang--Baxter property and singularity invariance]{Yang--Baxter property and singularity invariance}
\label{sec:YB}

To formulate a criterium for a quadrirational map to satisfy the
Yang--Baxter relation, we start with the parametrization of the
quadrirational maps in terms of their singularity sets following
Adler--Bobenko--Suris \cite{ABS1}. For simplicity we restrict our
considerations to the generic case corresponding to the
$({\mathcal M}\ddot{o}b)^4$-orbits of the maps of type $F_\I$. In
that case the correspon\-ding quadrirational map $F$ has the
singular set $\Sigma(F)=\{P_1,\ldots,P_4\}$ consisting of four
distinct points $P_i = (x_i, y_i) \in \mathbb{CP}^1\times
\mathbb{CP}^1$ (where the numerators and the denominators of the
both fractions in (\ref{maps}) simultaneously vanish).
Analogously, the map $F^{-1}$ has the singular set
$\Sigma(F^{-1})=\{Q_1,\ldots,Q_4\}$ consisting of four distinct
points $Q_i = (u_i, v_i) \in \mathbb{CP}^1\times \mathbb{CP}^1$.
Moreover, each of the companion maps $\bar R$ and $\bar R^{-1}$
has four distinct singular points, given by $(x_i,v_i)$ and
$(u_i,y_i)$, respectively. It can be shown that $F$ blows up the
points $P_i$ and blows down four exceptional curves~$C_i$ of
bidegree $(1,1)$ to the points $Q_i$. The exceptional curves can
be ordered in such a~way that $C_i$ passes through three points~$P_j$, $j\neq i$ (see Theorems~21 and~22 in~\cite{ABS1}).

Note that an ordered set $X=(x_1,\ldots,x_4)$ of 4 distinct points
in $\mathbb{CP}^1$ is in a one-to-one correspondence with the set
of pairs $(\a,\sigma)$, $\a\in\mathbb{C}$, $\sigma\in {\mathcal
M}\ddot{o}b$. Indeed, $\a$ is def\/ined as the cross ratio
\[
\alpha = q(X)=\frac{(x_1-x_2)(x_3-x_4)}{(x_2-x_3)(x_4-x_1)},
\]
and $\sigma$ is uniquely def\/ined by the condition that
\begin{gather*}
\sigma(x_1)=\infty, \qquad \sigma(x_2)=1,\qquad \sigma(x_3) = 0,\qquad
\sigma(x_4) = \a,
\end{gather*}
giving
\[
\sigma(x) = \frac{(x_1-x_2)(x-x_3)}{(x_3-x_2)(x-x_1)}.
\]
Let $X=(x_1,\dots, x_4)$, $Y = (y_1, \dots, y_4)$, $U = (u_1,
\dots, u_4)$, $V=(v_1, \dots, v_4)$ be four ordered quadruples of
elements of $\mathbb{CP}^1$. According to Theorem 18 in
\cite{ABS1}, they correspond to some quadrirational map if and
only if $U = \sigma(X)$, $ V=\tau (Y)$ for some $\sigma, \tau \in
{\mathcal M}\ddot{o}b$, that is, if $q(X)=q(U)$ and $q(Y)=q(V)$.
We denote the corresponding map by $F(\lambda, \mu)$, where
$\lambda=(X,U)$, $\mu=(Y,V)$ belong to the {\it admissible} set
\[
\Lambda=\big\{(X,U)\in(\mathbb{CP}^1)^4 \times (\mathbb{CP}^1)^4:\
U=\sigma(X)\ {\rm with}\ \sigma\in {\mathcal M}\ddot{o}b\big\}.
\]
Clearly, a simultaneous permutation of all four quadruples
$X$, $Y$, $U$, $V$ leads to the same map~$F$.

We now turn to the Yang--Baxter relation for the class of
quadrirational maps. In the formulation of our main result, we use
the following notation. Let $\pi \in S_4$ be a permutation of
indices $(1,\ldots,4)$. Then
\[
\Lambda_\pi=\big\{(X,U)\in\Lambda:\ U=X^{\pi}\Leftrightarrow\
u_i=x_{\pi(i)}, \ i=1,\dots, 4\big\}.
\]
Since in this def\/inition $U$ must be ${\mathcal
M}\ddot{o}b$-equivalent to $X$, the permutation $\pi$ must
preserve the cross ratio and thus belongs to the {\it Klein
subgroup} $K \subset S_4$ consisting of the identity $\rm Id$ and the
three commuting involutions
\[
\rho_1=(12)(34),\qquad \rho_2=(14)(23), \qquad \rho_3=(13)(24).
\]
Note that any two of $\rho_i$ are conjugated in $S_4$.

\begin{theorem}\label{main th}
The map $F(\lambda,\mu)$ satisfies the Yang--Baxter relation
\eqref{eq:YBcom} if and only if all parameters $\lambda_1$,
$\lambda_2$, $\lambda_3$ belong to the same subset $\Lambda_\pi
\subset \Lambda$, $\pi\in K$.
\end{theorem}

Note that for $F(\lambda,\mu)$ with $\lambda,\mu\in\Lambda_\pi$
the singular set of the inverse map consists of
$(u_i,v_i)=(x_{\pi_i},y_{\pi_i})$, and coincides with the singular
set of $F$, which can be considered as the following version of the singularity conf\/inement condition.

Let us say that a quadrirational map $F$ satisf\/ies {\it singularity invariance}\footnote{A closely related property called by Duistermaat~\cite{D} {\it geometric confinement}  was introduced in \cite{V}.}
property if the singular set of $F^{-1}$ coincides with that of the map $F.$
We see that the condition
$\lambda,\mu\in\Lambda_\pi$ is precisely the singularity invariance
 for the corresponding map $F(\lambda, \mu)$.

Thus, our main observation is that the Yang--Baxter property for
the quadrirational maps is equivalent to the singularity
invariance.  In particular, this explains why the canonical maps
$F_\I$--$F_\V$ satisfy the Yang--Baxter relation. For example, for
the map $F_\I$ we have:
\[
X=U=(\infty,1,0,\a),\qquad Y=V=(\infty,1,0,\beta),
\]
so that the singularity set is
\[
\Sigma(F_\I) =\{(\infty, \infty), (1,1), (0,0), (\a,
\b)\}=\Sigma\big(F_\I^{-1}\big).
\]
A quadrirational map $R$ with conf\/ined singularities determines a
permutation $\pi(R)$ of the singularity set $\Sigma(R)$. In the
$F_\I$ case this permutation is identity, but this need not be the
case in general.

For example, if we change in $F_\I$ the variables $y\mapsto\beta
y^{-1}$, $u \mapsto \alpha u^{-1}$, while keeping~$x$,~$v$ unchanged,
we obtain a non-equivalent map
\begin{gather}
\label{HI} \tag{$H_\I$}
 u =  y Q^{-1},\qquad v =  x Q, \qquad
 Q = \frac{(1-\b)xy+(\b-\a)y+\b(\a-1)} {(1-\a)xy+(\a-\b)x+\a(\b-1)},
\end{gather}
which we denote by $H_\I$. For this map we have:
\[
X=(\infty,1,0,\alpha),\qquad Y=(0,\beta,\infty,1), \qquad
U=(0,\alpha,\infty,1),\qquad V=(\infty,1,0,\beta).
\]
The singularity set is
\[
\Sigma(H_\I) =\{(\infty,0), (1,\b), (0,\infty), (\a,1)\} =
\Sigma\big(H_\I^{-1}\big),
\]
and the corresponding permutation is $\pi =\rho_3= (13)(24)$.
Thus, the map $H_\I$ has the singularity invariance property and,
correspondingly, it satisf\/ies the Yang--Baxter relation. It is
conjugated to the so-called Harrison map found in \cite{PTV}, see
also \cite{Harrison,TTX}. Similar maps can be constructed
for other canonical maps except $F_{\IV}$ (see the full list in
Section~\ref{section4} below). In the following we elucidate these issues in
detail.

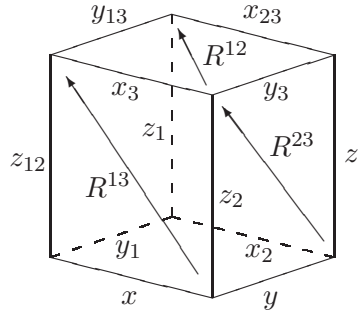
\begin{figure}[htbp]\setlength{\unitlength}{0.07em}
\centering
  \begin{picture}(200,170)(-80,-50)
   \dashline{5}(-60,-20)(0,0)(0,100)\dashline{5}(0,0)(80,-20)
   \path(-60,80)(-60,-20)(20,-40)(80,-20)(80,80)(0,100)
        (-60,80)(20,60)(80,80) \path(20,60)(20,-40)
   \put(75,-12){\vector(-3,4){50}}
   \put(13,-28){\vector(-2,3){65}}
   \put(17, 65){\vector(-1,2){14}}
   \put(-25,-43){$x$}   \put( 36,-20){$x_2$}
   \put(-30, 60){$x_3$} \put( 35, 98){$x_{23}$}
   \put( 45,-43){$y$}   \put(-28,-18){$y_1$}
   \put( 45, 60){$y_3$} \put(-40, 98){$y_{13}$}
   \put(85,27){$z$}     \put(23,7){$z_2$}
   \put(-15,40){$z_1$}  \put(-80,25){$z_{12}$}
   \put(-43,10){$R^{13}$}\put(15,75){$R^{12}$}\put(48,30){$R^{23}$}
 \end{picture}
  \caption{Right-hand side of Yang--Baxter equation; the left hand side is
  the composition of the maps corresponding to the other three faces of the cube.}
  \label{fig:YB}
\end{figure}

To prove the theorem recall that  according to Theorem 19 in
\cite{ABS1}, multidimensionally consistent systems of
quadrirational maps are characterized by matching of singularities
in the following way. (We restrict ourselves here again with the
generic case of the systems consisting of maps from the $(\mathcal
M\ddot{o}b)^4$-orbit of $F_\I$.) Let each edge of the
three-dimensional cube carry an ordered quadruple of elements of
$\mathbb{CP}^1$, such that the cross-ratios of the quadruples on
the opposite edges of each face coincide. We denote these
quadruples by the the same letters as the corresponding f\/ields on
Fig.~\ref{fig:YB}, but use the capital case. Attach to each face
of the cube a quadrirational map according to this data. Then this
system of maps is three-dimensionally consistent.

In order for the maps on the bottom and on the top faces to
coincide it is necessary and suf\/f\/icient that the corresponding
eight quadruples be related by
\[
(X_3,Y_3,X_{23},Y_{13})=(X,Y,X_2,Y_1)^{\pi_3}
\]
with some permutation $\pi_3\in S_4$. (This is to be read as the
simultaneous permutation of all four quadruples, so that
$X_3=X^{\pi_3}$ etc.) Equal cross-ratios for opposite edges
implies, as discussed above, that $\pi_3\in K$; in particular, it
is an involution. Similarly, each of the other two pairs consists
of the same maps if
\[
(Y_1,Z_1,Y_{13},Z_{12})=(Y,Z,Y_3,Z_2)^{\pi_1},\qquad
(X,Z_2,X_3,Z_{12})=(X_2,Z,X_{23},Z_1)^{\pi_2}
\]
with some involutive permutations $\pi_1,\pi_2\in K\subset S_4$.
These relations are equivalent to the following ones: f\/irst,
\begin{alignat*}{3}
&(X,X_2)\in\Lambda_{\pi_2},\qquad &&  (Y,Y_1)\in\Lambda_{\pi_1}, & \\
& (Y,Y_3)\in\Lambda_{\pi_3},\qquad && (Z,Z_2)\in\Lambda_{\pi_2}, & \\
& (X,X_3)\in\Lambda_{\pi_3},\qquad && (Z,Z_1)\in\Lambda_{\pi_1},&
\end{alignat*}
and second, that any two of the permutations $\pi_1$, $\pi_2$, $\pi_3$
commute. In order for the corresponding maps to be identif\/ied with
$F^{12}(\lambda_1,\lambda_2)$, $F^{23}(\lambda_2,\lambda_3)$, and
$F^{13}(\lambda_1,\lambda_3)$, it is required that
$\lambda_1=(X,X_2)=(X,X_3)$, $\lambda_2=(Y,Y_1)=(Y,Y_3)$, and
$\lambda_3=(Z,Z_1)=(Z,Z_2)$. Therefore,  we arrive at the
condition $\pi_1=\pi_2=\pi_3$, which proves the theorem.

When $\pi= (13)(24)$ we arrive at the map $H_\I$ def\/ined above.
For other nontrivial elements of the Klein group we have dif\/ferent families of maps,
which however turn out to be {\it equivalent as families} if we allow also change of parameters
(see below).  This corresponds to the fact that the nontrivial elements of the Klein
subgroup are conjugated to one another.

Thus, in the generic case we have two non-equivalent 8-parameter families of quadrirational Yang--Baxter maps:
\begin{alignat*}{3}
& R(X, Y) =  F((X,X), (Y,Y))  = \sigma^{-1}\times\tau^{-1} F_\I(\a, \b)  \sigma\times\tau
 \qquad && (F_\I\mbox{-family}), & \\
& R(X, Y) =  F((X,X^{\rho}), (Y,Y^{\rho}))= \sigma^{-1}\times\tau^{-1}H_\I (\a, \b) \sigma\times\tau   \qquad && (H_\I\mbox{-family}),&
\end{alignat*}
where $\rho= (13)(24)$, $X, Y \in (\mathbb{CP}^1)^4$ are ordered sets of 4 distinct points in $\mathbb{CP}^1$
and $\sigma, \tau \in {\mathcal M}\ddot{o}b.$

Of course these maps are equivalent, in the sense of Proposition~\ref{Pr}, to the maps $F_\I$ and $H_\I$ respectively. Next, we are going to show that  additional families correspond to the symmetries of the pencils of quadrics.

\section[Geometric symmetries and additional families of Yang-Baxter maps]{Geometric symmetries and additional families\\ of Yang--Baxter maps}
\label{section4}

The relation between the maps $F_\I$ and $H_\I$ can be also
established with the help of the following statement. We say that
a family of involutions $\sigma(\lambda)$ of the set $\mathbb X$
is a {\it symmetry} of a reversible YB map $R(\lambda,\mu):\mathbb
X\times\mathbb X\to\mathbb X\times\mathbb X$ if
\begin{gather*}
\sigma (\lambda) \times \sigma (\mu)   R(\lambda,\mu)= R(\lambda,\mu)
\sigma (\lambda) \times \sigma (\mu).
\end{gather*}
The following simple but important statement is easily proved by a
direct check.

\begin{proposition}\label{lemma}
If $\sigma(\lambda)$ is a symmetry of a reversible YB map $R(\lambda,\mu)$,
then
\begin{gather*}
R^{\sigma} = \sigma (\lambda)\times {\rm Id} \, R(\lambda,\mu)\, {\rm Id} \times \sigma (\mu)
\end{gather*}
is also a reversible YB map.
\end{proposition}

\begin{figure}[htbp]
\centerline{\includegraphics[width=7cm]{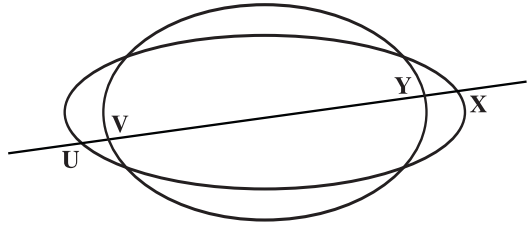}}
\caption{A quadrirational map on a pair of conics.}\label{fig:XYUV}
\end{figure}

A nice way to f\/ind symmetries of the canonical maps $F_\I$--$F_\V$ is
provided by their geometric interpretation given in \cite{ABS1}
and illustrated on Fig.~\ref{fig:XYUV}. Consider two conics $Q_1$,
$Q_2$ on the plane and let $\mathbf{X}\in Q_1$, $\mathbf{Y}\in
Q_2$. Set $\mathbf{U}$ ($\mathbf{V}$) to be the second
intersection point of the line $\overline{\mathbf{X}\mathbf{Y}}$
with $Q_1$ (resp.\ with $Q_2$). Then the maps
$F_\I$--$F_\V:(x,y)\mapsto(u,v)$ are just the coordinate
representations of the map ${\mathcal
R}:(\mathbf{X},\mathbf{Y})\mapsto(\mathbf{U},\mathbf{V})$ in some
special rational parametrization of the conics, $\a$~and~$\b$
being the cross-ratios of the corresponding four intersection
points. The f\/ive maps $F_\I$--$F_\V$ correspond to the f\/ive possible
types $\I$--$\V$ of intersection of two conics:
\begin{itemize} \itemsep=0pt
\item[\I:] four simple intersection points; \item[\II:] two simple
intersection points and one point of tangency; \item[\III:] two
points of tangency; \item[\IV:] one simple intersection point and
one point of the second order tangency; \item[\V:] one point of
the third order tangency.
\end{itemize}

\begin{figure}[htbp]
\centerline{\includegraphics[width=7cm]{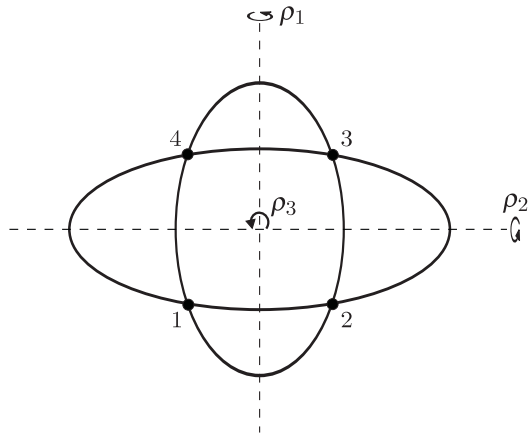}}
\caption{Involutions of the pair of quadrics in the generic case.}
\label{fig:klein}
\end{figure}

Let us examine f\/irst the symmetries of the $F_\I$ map. The pair of
quadrics in the $w_1-w_2$ plane, in generic position, can be
reduced by a suitable projective transformation to
\[
Q_1:\ w_2(w_2-1)=\alpha w_1(w_1-1),\qquad Q_2:\ w_2(w_2-1)=\beta
w_1(w_1-1).
\]
It has an obvious symmetry group $K=\mathbb Z_2\times\mathbb
Z_2=\{{\rm Id},\rho_1,\rho_2,\rho_3\}$, where
\begin{gather*}
\rho_1(w_1,w_2)=(1-w_1,w_2),\!\!\qquad \rho_2(w_1,w_2)=(w_1,1-w_2),
\!\!\qquad \rho_3(w_1,w_2)=(1-w_1,1-w_2)\!
\end{gather*}
(see Fig.~\ref{fig:klein}).

\begin{figure}[htbp]
\centerline{\includegraphics[width=7cm]{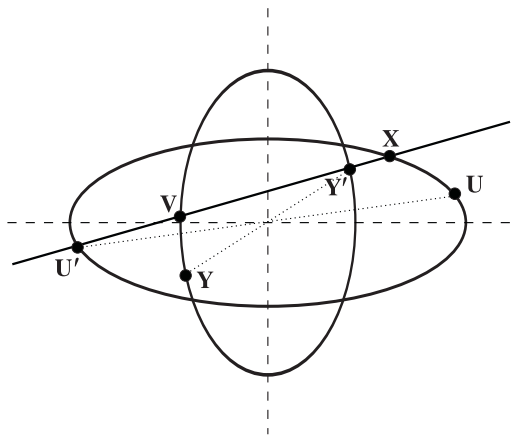}}
\caption{Geometric interpretation of $H_I$ map.}
\label{fig:HI}
\end{figure}

In the parametrization
\[
w_1=\frac{x-\a}{x^2-\a},\qquad w_2=\frac{x(x-\a)}{x^2-\a}
\]
for the f\/irst quadric and similarly for the second one with $\a$
replaced by $\b$, the symmetry $\rho_3$ corresponds to
\[
\sigma(\a): \ x \mapsto \frac{\a}{x}.
\]
By applying Proposition~\ref{lemma} to the map $F_\I$ with this
symmetry, we obtain the map $H_\I$ (see Fig.~\ref{fig:HI}).

Application of Proposition~\ref{lemma} to $F_\I$ with the symmetry $\rho_1$, which
corresponds to
\[
\tau(\a): \ x \mapsto \frac{x-\alpha}{x-1},
\]
leads to the YB map
\begin{gather}\label{HI2}
u=\frac{\a y(\b -x-y+xy)}{\a\b - \b x -\a y + \b xy}, \qquad
v=\frac{\b x(\a-x-y+xy)}{\a\b - \b x -\a y+ \a xy}.
\end{gather}
It is, however, equivalent to the map $H_\I$ via
the conjugation by $x\mapsto(x-1)/x$ and the reparametrization
$\a\mapsto(\a-1)/\a$. Similarly, the third symmetry $\rho_2$ leads
to the family which is related to the map (\ref{HI2}) by
conjugation $x\mapsto 1/x$ and reparametrization $\a\mapsto 1/\a$.

In the degenerate cases $F_{\II}$ and $F_\V$ the symmetry group
breaks down to $\mathbb Z_2$, so we have one additional family for
each case. In case $F_{\III}$ we can choose the equations of
quadrics as
\[
Q_1:\ \alpha w_1(w_1-1)=w_2^2, \qquad Q_2:\ \beta
w_1(w_1-1)=w_2^2,
\]
and we still have $\mathbb Z_2 \times \mathbb Z_2$ symmetry group.
One can check that the symmetries
$(w_1,w_2)\mapsto(1-w_1,\pm w_2)$  lead to equivalent families, so
we have two new families in this case, depending whether the points of tangency are interchanged ($B$-family) or not ($A$-family). In the case $F_{\IV}$ we
have no non-trivial symmetries. This leads to the following list.

\begin{theorem}
Up to equivalence \eqref{equiva} and reparametrization there are $10$ families of
quadrirational reversible Yang--Baxter maps of subclass {\rm [2:2]}: the $5$~families
$F_\I$--$F_\V$ from {\rm \cite{ABS1}} and the following $5$ additional
families
\begin{alignat}{4}
\label{H1}\tag{$H_\I$}
 & u= yQ^{-1},\qquad &&
  v= x Q,\qquad &&
  Q= \frac{(1-\b)xy+(\b-\a)y+\b(\a-1)} {(1-\a)xy+(\a-\b)x+\a(\b-1)},&\\
\label{H2}\tag{$H_\II$}
&  u= y Q^{-1},\qquad &&
  v= x Q,\qquad &&
  Q=\frac{\a + (\b-\a) y - \b xy}{\b + (\a-\b)x -\a xy},&\\
\label{H3A}\tag{$H_\III^A$}
&  u= \frac{y}{\a} Q,\qquad &&
  v= \frac{x}{\b} Q,\qquad &&
 Q= \frac{\a x + \b y}{x+y},&\\
 \label{H3B}\tag{$H_\III^B$}
&  u= y Q^{-1},\qquad &&
  v=  x Q,\qquad &&
 Q= \frac{\a xy + 1}{\b xy+1},&\\
\label{H5}\tag{$H_\V$}
&  u= y-P,\qquad &&
  v= x+P,\qquad &&
  P= \frac{\a-\b}{x+y}.&
\end{alignat}
\end{theorem}

The map $H_\V$ is the well-known Adler map (see \cite{V1}). The
map $H_{\III}^B$ appeared in~\cite{PTV} in relation with the
discrete potential KdV equation and in~\cite{KNW} in relation with
the discrete KP hierarchy.

Note that $H_{\III}^A$ and $H_\V$ can be obtained from $F_{\III}$
and $F_\V$, respectively, by the change $u\mapsto -u$, $y \mapsto
-y$.

The maps $H_{\I}$, $H_{\II}$ have convenient subtraction-free representatives:
\begin{alignat}{3}
\label{H1plus}\tag{$H_{\I}^+$}
&  u= \frac{y}{A} \frac{B+A x+ B y +ABxy}{1+x+y+Bxy},\qquad &&
  v= \frac{x}{B} \frac{A+A x+ B y +ABxy}{1+x+y+Axy},&\\
\label{H2plus}\tag{$H_{\II}^+$}
&  u= \frac{y}{\a} \frac{\a x+ \b y +\b}{x+y+1},\qquad &&
  v= \frac{x}{\b} \frac{\a x+ \b y +\a}{x+y+1},&
 \end{alignat}
where $A=1-\a$, $B=1-\b$. They are related to the maps $H_{\I}$,
$H_{\II}$ above by conjugation viz.\ $x\mapsto 1/(x-1)$ and
similarly for all other variables. Note that the maps $H_{\III}^A$
and $H_{\III}^B$ are subtraction-free also, so $H_{\I}^+$, $H_{\II}^+$, $H_{\III}^A$ and $H_{\III}^B$  have a natural
``tropical'' (or ultra-discrete) limit (cf.~\cite{KNW}). Moreover
there is an obvious (singular) limit procedure starting from
$H_{\I}^+$ in order to obtain $H_{\II}^+$ and $H_{\III}^A$.

Although the families of maps $F$ and $H$ are not equivalent it is
not clear if they lead to dif\/ferent  transfer dynamics~\cite{V1}.
For example, both $F_\V$ and $H_\V$ are related to the discrete
(potential) KdV equation, but via dif\/ferent symmetries of it (see
\cite{PTV}). However, we know that considering other symmetries of
discrete KdV equation one can obtain the families $F_{\IV}$ and
even~$H_\I$ (see~\cite{PTV}).

Note also that although all the maps $F_\I$--$F_\V$, $H_\I$--$H_\V$ are the involutions, the corresponding transfer-dynamics \cite{V1} is non-trivial and deserves further investigation.

\subsection*{Acknowledgements}
This work had been started in May 2007, when one of the authors
(APV) visited Patras within the Erasmus-Socrates exchange programme,
and completed at the Isaac Newton Institute
for Mathematical Sciences in Cambridge during the programme on
Discrete Integrable Systems in the spring semester 2009.
The work of APV was also partially
supported by the European RTN ENIGMA (contract MRTN-CT-2004-5652)
and EPSRC (grant EP/E004008/1). We are grateful to V.~Adler for
useful comments.

\pdfbookmark[1]{References}{ref}
\LastPageEnding

\end{document}